\documentclass{amsproc}

\usepackage{amssymb}

\usepackage{graphicx}

\usepackage[cmtip,all]{xy}

\usepackage[all]{xy}
\usepackage{amsthm}
\usepackage{amscd}
\usepackage{amsmath}
\usepackage{amsfonts}
\usepackage{amssymb}
\usepackage{stmaryrd}
\usepackage{graphicx}%
\usepackage{hyperref}
\hypersetup{
    bookmarks=true,         
    unicode=false,          
    pdftoolbar=true,        
    pdfmenubar=true,        
    pdffitwindow=false,     
    pdfstartview={FitH},    
    pdftitle={My title},    
    pdfauthor={Author},     
    pdfsubject={Subject},   
    pdfcreator={Creator},   
    pdfproducer={Producer}, 
    pdfkeywords={keyword1, key2, key3}, 
    pdfnewwindow=true,      
    colorlinks=false,       
    linkcolor=red,          
    citecolor=green,        
    filecolor=magenta,      
    urlcolor=cyan           
}

\newtheorem{theorem}{Theorem}[section]
\newtheorem{lemma}[theorem]{Lemma}

\theoremstyle{definition}
\newtheorem{definition}[theorem]{Definition}
\newtheorem{example}[theorem]{Example}
\newtheorem{proposition}[theorem]{Proposition}

\theoremstyle{remark}
\newtheorem{remark}[theorem]{Remark}

\numberwithin{equation}{section}

\begin{document}

\title{Quasi-theories and their equivariant orthogonal spectra}



\author{Zhen Huan}

\address{Zhen Huan, Department of Mathematics,
Sun Yat-sen University, Guangzhou, 510275 China} \curraddr{}
\email{huanzhen84@yahoo.com}


\subjclass[2010]{Primary 55}

\date{}

\begin{abstract}
In this paper we construct orthogonal $G-$spectra up to a weak equivalence for the quasi-theory $QE_{n, G}^*(-)$ corresponding to certain cohomology theories $E$. The construction of the orthogonal $G-$spectrum for quasi-elliptic cohomology can be applied to the constructions for quasi-theories. 
\end{abstract}

\maketitle


\section{Introduction}

In \cite{HuanSpec} we construct a functor $\mathcal{Q}$ from the category of
orthogonal ring spectra to the category of $\mathcal{I}_G-$FSP. If $E$ is a global cohomology theory, $\mathcal{Q}(E)$
weakly represents the
cohomology theory \begin{equation}QE^*_G(X):=\prod_{\sigma\in
G^{tors}_{conj}}E^*_{\Lambda(\sigma)}(X^{\sigma})=\bigg(\prod_{\sigma\in
G^{tors}}E^*_{\Lambda(\sigma)}(X^{\sigma})\bigg)^G.\label{defintroE2}\end{equation}
The image of global $K-$spectrum is a $\mathcal{I}_G-$FSP representing quasi-elliptic cohomology up to a weak equivalence.

Quasi-elliptic cohomology is a variant of
elliptic cohomology theories, which is the generalized elliptic cohomology theory associated to the Tate curve $Tate(q)$ over Spec$\mathbb{Z}((q))$  [Section 2.6,
\cite{AHS}].  Quasi-elliptic cohomology is defined over Spec$\mathbb{Z}[q^{\pm}]$.
Inverting $q$ allows us to define a sufficiently non-naive equivariant cohomology theory and to interpret some constructions more easily.
Its relation with Tate K-theory is \begin{equation}QEll^*_G(X)\otimes_{\mathbb{Z}[q^{\pm}]}\mathbb{Z}((q))=(K^*_{Tate})_G(X)
\label{tateqellequiv}\end{equation}

Motivated by quasi-elliptic cohomology, we construct quasi-theories $QE_{n, G}^*(-)$ in \cite{Huanquasi}. Quasi-elliptic cohomology, the theories $QE^*_G(-)$ defined in (\ref{defintroE2}) and the generalized quasi-elliptic cohomology in Example \ref{generalizedquasi} are all special cases of quasi-theories.

In this paper we show that the idea of constructing the functor $\mathcal{Q}$ can be applied to construct a family of functors $\mathcal{Q}_n$ from the category of orthogonal ring spectra to the category of $\mathcal{I}_G-$FSP. Especially, the functor $\mathcal{Q}_1$ is $\mathcal{Q}$.
In other words, we construct a $\mathcal{I}_G-$FSP representing $QE_{n, G}^*(-)$ up to weak equivalence for each positive integer $n$ and each compact Lie group $G$.

In this paper we show the construction of functors $\mathcal{Q}_n$. The idea is analogous to the construction of $\mathcal{Q}$ in \cite{HuanSpec}. For the readers' convenience, we still include all the details in this paper. In Section \ref{definequasi} we recall the definition and examples of quasi-theories. In Section \ref{newcatsp}
we recall a category of orthogonal $G-$spectra introduced in \cite{HuanSpec}. In Section \ref{equisp} we construct a space $QE_{G, n, m}$ representing the $m-$th
$G-$equivariant quasi-theory $QE^m_{n, G}(-)$ up to a weak equivalence. In Section \ref{GorthospectraQEll} we construct a $\mathcal{I}_G-$FSP representing $QE_{n, G}^*(-)$ up to weak equivalence for certain cohomology theories $E$ and construct the functors $\mathcal{Q}_n$.
In the appendix, we construct some faithful group representations needed in the construction of the $\mathcal{I}_G-$FSP.

\subsection{Acknowledgement}
I would like to thank my PhD advisor Charles Rezk. Under his direction I constructed equivariant orthogonal spectra for quasi-elliptic cohomology, which is a special case of quasi-theories. In addition, he suggested the project on quasi-theories to me. I would like to thank Matthew Ando for encouraging me to finish the projects.

\section{The Quasi-theory $QE_{n, G}^*(-)$}\label{definequasi}

In this section we recall the quasi-theories. The main reference for that is \cite{Huanquasi}.

Let $G$ be a compact Lie group and $n$ denote a positive integer. Let $G^{tors}_{conj}$ denote a set of representatives of
$G-$conjugacy classes in the set $G^{tors}$  of torsion elements in $G$. Let $G^{n}_{z}$ denote set
$$\{\sigma=(\sigma_1, \sigma_2, \cdots \sigma_n )| \sigma_i\in G^{tors}_{conj}, [\sigma_i, \sigma_j]\mbox{   is  the   identity  element in  }G\}.$$

Let $\sigma=(\sigma_1, \sigma_2, \cdots \sigma_n)\in G^n_z$.  Define
\begin{align}C_G(\sigma)&:=\bigcap\limits_{i=1}^nC_G(\sigma_i); \label{Csigmadef}\\ \Lambda_G(\sigma)&:= C_G(\sigma)\times
\mathbb{R}^n/\langle (\sigma_1, -e_1), (\sigma_2, -e_2), \cdots (\sigma_n, -e_n)\rangle.\label{lambdadef}\end{align}
where $C_G(\sigma_i)$ is the centralizer of each $\sigma_i$ in $G$ and $\{e_1, e_2, \cdots e_n\}$ is a basis of $\mathbb{R}^n$.
Let $q:\mathbb{T}\longrightarrow U(1)$ denote the representation $t\mapsto e^{2\pi i t}$. Let $q_i=1\otimes\cdots\otimes q\otimes\cdots\otimes 1: \mathbb{T}^{n}\longrightarrow U(1)$ denote the tensor product with $q$ at the $i-$th position and trivial representations at other position. The representation ring $$R(\mathbb{T}^{n})\cong R(\mathbb{T})^{\otimes n}=\mathbb{Z}[q_1^{\pm}, \cdots q_n^{\pm}].$$

We have the exact sequence  \begin{equation}1\longrightarrow C_G(\sigma)\longrightarrow \Lambda_G(\sigma)\buildrel{\pi}\over\longrightarrow \mathbb{T}^{n}\longrightarrow 0 \end{equation}
where the first map is $g\mapsto [g, 0]$  and the second map is $\pi([g, t_1, \cdots t_n])=(e^{2\pi i t_1}, \cdots e^{2\pi it_n}).$
Then the map $\pi^*: R(\mathbb{T}^{n})\longrightarrow R\Lambda_G(\sigma)$ equips the representation ring $R\Lambda_G(\sigma)$ the structure as an
 $R(\mathbb{T}^{n})-$module.

This is Lemma 3.1 \cite{Huanquasi} presenting the relation between $RC_G(\sigma)$ and $R\Lambda_G(\sigma)$.

\begin{lemma} $\pi^*: R(\mathbb{T}^{n})\longrightarrow R\Lambda_G(\sigma)$ exhibits $R\Lambda_G(\sigma)$ as a free $R(\mathbb{T}^{n})-$module.

There is an $R(\mathbb{T}^{n})-$basis of $R\Lambda_G(\sigma)$
given by irreducible representations $\{V_{\lambda}\}$, such that
restriction $V_{\lambda}\mapsto V_{\lambda}|_{C_G(\sigma)}$ to $C_G(\sigma)$
defines a bijection between $\{V_{\lambda}\}$ and the set
$\{\lambda\}$ of irreducible representations of
$C_G(\sigma)$.\label{cl}\end{lemma}

\begin{definition}For equivariant cohomology theories $\{E_{H}^*\}_H$ and any $G-$space $X$, the corresponding quasi-theory $QE_{n, G}^*(X)$ is defined to be
$$\prod_{\sigma\in
G^{n}_{z}}E^*_{\Lambda_G(\sigma)}(X^{\sigma}).$$\label{qedef}\end{definition}

\begin{example}[Motivating example: Tate K-theory  and quasi-elliptic cohomology]
Tate $K-$theory is the generalized elliptic cohomology associated to the Tate curve. The elliptic cohomology theories form a sheaf of cohomology theories over the moduli stack of elliptic curves
$\mathcal{M}_{ell}$. Tate K-theory over Spec$\mathbb{Z}((q))$ is obtained when we restrict it to a punctured completed neighborhood of the cusp at $\infty$, i.e. the Tate
curve $Tate(q)$ over Spec$\mathbb{Z}((q))$  [Section 2.6,
\cite{AHS}].  The divisible group associated to Tate K-theory is $\mathbb{G}_m\oplus \mathbb{Q}/\mathbb{Z}$.
The relation between Tate K-theory and string theory is better
understood than most known elliptic cohomology theories. In addition, Tate K-theory has the closest ties to Witten's original insight that the elliptic cohomology of  a space $X$ is related to the
$\mathbb{T}-$equivariant K-theory of the free loop space
$LX=\mathbb{C}^{\infty}(S^1, X)$ with the circle $\mathbb{T}$
acting on $LX$ by rotating loops. Ganter gave a careful interpretation in Section 2, \cite{Gan07} of this statement that the
definition of $G-$equivariant Tate K-theory for finite groups $G$
is modelled on the loop space of a global quotient orbifold.

Other than the theory over Spec$\mathbb{Z}((q))$, we can define variants of Tate K-theory over Spec$\mathbb{Z}[q]$ and Spec$\mathbb{Z}[q^{\pm}]$ respectively. The theory over Spec$\mathbb{Z}[q^{\pm}]$
is of especial interest.
Inverting $q$ allows us to define a sufficiently non-naive equivariant cohomology theory and to interpret some constructions more easily in terms of extensions of groups over
the circle. The resulting cohomology theory is called quasi-elliptic cohomology \cite{Rez11}\cite{Huanthesis}\cite{Huansurvey}.
Its relation with Tate K-theory is \begin{equation}QEll^*_G(X)\otimes_{\mathbb{Z}[q^{\pm}]}\mathbb{Z}((q))=(K^*_{Tate})_G(X)
\label{tateqellequiv}\end{equation} which also reflects the geometric nature of the Tate curve.
$QEll^*_{\mathbb{T}}(\mbox{pt})$ has a direct interpretation in terms of the Katz-Mazur group scheme $T$ [Section 8.7, \cite{KM85}].
The idea of quasi-elliptic cohomology is  motivated by Ganter's
construction of Tate K-theory  \cite{Dev96}.
It is not an elliptic cohomology  but a more robust and algebraically simpler treatment of
Tate K-theory. This new theory can be interpreted in a neat form
by equivariant K-theories. Some
formulations in it can be generalized to equivariant cohomology
theories other than Tate K-theory.

Quasi-elliptic cohomology $QEll^*_G(-)$ is exactly the quasi-theory $QK_{1, G}^*(-)$ in Definition \ref{qedef}.
\end{example}

\begin{example}[Generalized Tate K-theory and generalized quasi-elliptic cohomology]

In Section 2 \cite{Gan07} Ganter  gave an interpretation of $G-$equivariant Tate K-theory for finite groups $G$
by the loop space of a global quotient orbifold. Apply the loop construction $n$ times, we can get the $n-$th generalized Tate K-theory. The divisible group associated to it is
$\mathbb{G}_m\oplus (\mathbb{Q}/\mathbb{Z})^n$.

With quasi-theories, we can get a neat expression of it. Consider the quasi-theory
$$QK_{n, G}^*(X)=\prod_{\sigma\in
G^{n}_{z}}K^*_{\Lambda_G(\sigma)}(X^{\sigma}).$$
$QK_{n, G}^*(X)\otimes_{\mathbb{Z}[q^{\pm}]^{\otimes n}}\mathbb{Z}((q))^{\otimes n}$ is isomorphic to the $n-$th generalized Tate K-theory.

\label{generalizedquasi}
\end{example}

\section{A new category of orthogonal $G-$spectra}\label{newcatsp}

It is difficult to construct a concrete representing spectrum for elliptic
cohomology.
In Section 4 \cite{HuanSpec} we formulate a new category of spectra with larger class of weak equivalence than that in \cite{MM}.
In Section 6 \cite{HuanSpec} we construct an orthogonal
$G-$spectrum for any compact Lie group $G$ representing $QEll^*_G(-)$ in this new category of
orthogonal $G-$spectra.

First we recall the category of orthogonal
$G-$spectra in \cite{MM} and the category $GwS$ that we will work in.
The weak equivalence of interest is the $\pi_*-$isomorphism.
\begin{definition} For subgroups $H$ of $G$ and integers $q$, define the homotopy groups $\pi^H_q(X)$ of a $G-$prepsectrum $X$ by
\begin{equation} \pi^H_q(X)= colim_V \pi^H_q(\Omega^V X(V))\mbox{  if  } q\geq 0,\end{equation} where $V$ runs over the indexing $G-$spaces in the chosen universe, and
\begin{equation} \pi^H_{-q}(X)= colim_{V\supset\mathbb{R}^q} \pi^H_0(\Omega^{V-\mathbb{R}^q} X(V))\mbox{  if  } q> 0.\end{equation}

A map $f: X\longrightarrow Y$ of $G-$prespectra is a $\pi_*-$isomorphism if it induces isomorphisms on all homotopy groups.

A map of orthogonal $G-$spectra is a $\pi_*-$isomorphism if its underlying map of $G-$prespectra is a $\pi_*-$isomorphism.  \label{pistar}
\end{definition}

\begin{definition}The category $GwS$ is the homotopy category of the category of orthogonal $G-$spectra with the weak equivalence defined by
\begin{equation}X\sim Y\mbox{  if   }\pi_0^H(X(V))=\pi_0^H(Y(V)),\label{weaksensp12}\end{equation} for each   faithful  $G-$representation  $V$ and any closed subgroup $H$ of $G$.

An orthogonal $G-$spectrum $X$ in $GwS$ is said to represent a theory $H^*_G$ if we have a natural map \begin{equation}\pi_0^H(X(V))=H^V_G(G/H),\label{weaksenspV}\end{equation} for each   faithful  $G-$representation  $V$ and any closed subgroup $H$ of $G$. \label{gwsdef}
\end{definition}

\begin{lemma} If a map $f:X\longrightarrow Y$ of orthogonal $G-$spectra induces isomorphisms (\ref{weaksensp12})  on the homotopy groups, i.e. \begin{equation}f: \pi_0^H(X(V))\buildrel\cong\over\longrightarrow \pi_0^H(Y(V))\label{wespdef}\end{equation} for each   faithful  $G-$representation  $V$,  and any closed subgroup $H$ of $G$, then $f$ is a $\pi_*-$isomorphism. \end{lemma}

We will work in the homotopy category of the category of orthogonal $G-$spectra with the weak equivalence defined in
(\ref{wespdef}). This homotopy category $Gw\mathcal{T}$ is smaller than the homotopy category of orthogonal $G-$spectra that we usually talked about, where the weak equivalence involved is the $\pi_*-$isomorphism. However, it seems the information that each object contains is enough to define an equivariant cohomology theory.


The homotopical adjunction below is a way to describe the relation between $G-$equivariant homotopy theory and those equivariant homotopy theory for its closed subgroups.
It is introduced in Definition 4.4 \cite{HuanSpec}.

\begin{definition}[homotopical adjunction] Let $H$ and $G$ be two compact Lie groups.
Let \begin{equation} L: G\mathcal{T} \longrightarrow H\mathcal{T}
\mbox{ and } R: H\mathcal{T} \longrightarrow
G\mathcal{T}\end{equation} be two functors between the category of $G-$spaces and that of $H-$spaces. A
\textit{left-to-right homotopical adjunction} is a natural map
\begin{equation}\mbox{Map}_H(LX,Y)\longrightarrow
\mbox{Map}_G(X,RY),\end{equation} which is a weak equivalence of
spaces when $X$ is a $G-$CW complex.

Analogously, a \textit{right-to-left homotopical adjunction} is a
natural map
\begin{equation}
\mbox{Map}_G(X,RY)\longrightarrow \mbox{Map}_H(LX,Y)\end{equation}
which is a weak equivalence of spaces when $X$ is a $G-$CW
complex.

$L$ is called a \textit{homotopical left adjoint} and $R$ a
\textit{homotopical right adjoint}. \label{holr}
\end{definition}

\section{Equivariant spectra}\label{equisp}

In this section,
we construct a space $QE_{G, n, m}$ representing the $m-$th
$G-$equivariant quasi-theory $QE^m_{n, G}(-)$ up to a weak equivalence.

Let $G$ be a compact Lie group and $\sigma\in G^{n}_z$.
Let $\Gamma$ denote the subgroup $\langle\sigma_1, \cdots \sigma_n\rangle$ of $G$. Let $$S_{G, \sigma}:=\mbox{Map}_{\Gamma}(G, \ast_{K}E(\Gamma/K))$$ where
$\ast$ denotes the join,  $K$ goes over all the maximal subgroups of
$\Gamma$ and $E(\Gamma/K)$ is the universal
space of the abelian group $\Gamma/K$.

\begin{lemma}

For any closed subgroup $H\leqslant G$, $S_{G, \sigma}$ satisfies
\begin{equation}S_{G, \sigma}^H\simeq\begin{cases}\mbox{pt},&\text{if for any $b\in G$,  $b^{-1}\Gamma b\nleq
H$;}\\
\emptyset,&\text{if there exists $b\in G$ such that
$b^{-1}\Gamma b \leqslant
H$.}\end{cases}\end{equation}\label{sf}
\end{lemma}

\begin{proof}
\begin{equation}S_{G,\sigma}^H=\mbox{Map}_{\Gamma}(G/H, \ast_{K}E(\Gamma/K)).\end{equation}

If there exists an $b\in G$ such that $b^{-1}\Gamma b
\leqslant H$, it is equivalent to that there exists points in
$G/H$ that can be fixed by $\Gamma$. But there are no points in
$\ast_{K}E(\Gamma/K)$ that can be fixed by the whole group $\Gamma$. So there
is no $\Gamma-$equivariant map from $G/H$ to
$\ast_{K}E(\Gamma/K)$. In this case $S_{G,\sigma}^H$ is
empty.

If for any $b\in G$, $b^{-1}\Gamma b\nleq H$, it is
equivalent to say that there are no points in $G/H$ that can be
fixed by $\Gamma$. Any proper subgroup $L$ of $\Gamma$ is contained in some maximal subgroup of $\Gamma$. $(\ast_{K}E(\Gamma/K))^{L}$ is the join of several
contractible spaces $E(\Gamma/K)^{L}$.
Thus, it is contractible. So all the homotopy groups
$\pi_n((\ast_{K}E(\Gamma/K))^{L})$ are
trivial. For any $n\geq 1$ and any $L-$equivariant
map $$f: (G/H)^n\longrightarrow \ast_{K}E(\Gamma/K)$$from the $n-$skeleton of $G/H$, the obstruction
cocycle is zero.

Then by equivariant obstruction theory, $f$ can be extended to the
$(n+1)-$cells of $G/H$, and any two extensions $f$ and $f'$ are
$\Gamma-$homotopic.


So in this case  $S_{G,\sigma}^H$ is contractible.
\end{proof}

\begin{theorem}
A homotopical
right adjoint of the functor $L_{\sigma}: G\mathcal{T}\longrightarrow
C_G(\sigma)\mathcal{T},\mbox{ } X\mapsto X^{\sigma}$ from the category of $G-$spaces to that of $C_G(\sigma)-$spaces is
\begin{equation}R_g:
C_G(\sigma)\mathcal{T}\longrightarrow G\mathcal{T}\mbox{,   }
Y\mapsto\mbox{Map}_{C_G(\sigma)}(G, Y\ast S_{C_G(\sigma), \sigma}).\end{equation}
\label{main}\end{theorem}
\begin{proof}
Let $H$ be any closed subgroup of $G$.

First we show given a $C_G(\sigma)-$equivariant map $f:
(G/H)^{\sigma}\longrightarrow Y$, it extends uniquely up to
$C_G(\sigma)-$homotopy to a $C_G(\sigma)-$equivariant map $$\widetilde{f}:
G/H\longrightarrow Y\ast S_{C_G(\sigma), \sigma}.$$ $f$ can be viewed as a
map $(G/H)^{\sigma}\longrightarrow Y\ast S_{C_G(\sigma), \sigma}$ by composing with
the inclusion of one end of the join
$$Y\longrightarrow Y\ast S_{C_G(\sigma), \sigma},\mbox{           } y\mapsto (1y, 0).$$

If $bH\in (G/H)^{\sigma}$,    define $\widetilde{f}(bH): =f(bH).$

If $bH$ is not in $(G/H)^{\sigma}$, its stabilizer group does not contain
$\Gamma$.  By Lemma \ref{sf}, for any subgroup $L$ of its stabilizer group,  $S_{C_G(\sigma),
\sigma}^L$ is contractible. So $(Y\ast S_{C_G(\sigma), \sigma})^{L}=Y^{L}\ast
S_{C_G(\sigma), \sigma}^{L}$ is contractible. In other words, if $L$ occurs
as the isotropy subgroup of a point outside $(G/H)^{\sigma}$,
$\pi_n((Y\ast S_{C_G(\sigma), \sigma})^{L})$ is trivial. By equivariant
obstruction theory, $f$ can extend
to a $C_G(\sigma)-$equivariant map $\widetilde{f}:
G/H\longrightarrow Y\ast S_{C_G(\sigma), \sigma}$, and any two extensions
are $C_G(\sigma)-$homotopy equivalent. In addition, $S_{C_G(\sigma), \sigma}^{\sigma}$
is empty. So the image of the restriction of any map
$G/H\longrightarrow Y\ast S_{C_G(\sigma),\sigma}$ to the subspace $(G/H)^{\sigma}$
is contained in the end $Y$ of the join.

Thus, $\mbox{Map}_{C_G(\sigma)}((G/H)^{\sigma}, Y)$ is weak equivalent to
$\mbox{Map}_{C_G(\sigma)}(G/H, Y\ast S_{C_G(\sigma), \sigma})$.

Moreover, we have the equivalence by adjunciton
\begin{equation}
\mbox{Map}_{G}(G/H, \mbox{           }\mbox{Map}_{C_G(\sigma)}(G, Y\ast
S_{C_G(\sigma), \sigma}))\cong \mbox{Map}_{C_G(\sigma)}(G/H, \mbox{ } Y\ast
S_{C_G(\sigma), \sigma})\end{equation} So we get
\begin{equation}R_{\sigma}Y^H=\mbox{Map}_{G}(G/H, R_{\sigma}Y)\backsimeq
\mbox{Map}_{C_G(\sigma)}((G/H)^{\sigma},
Y)\label{weakorbitbasic}\end{equation} Let $X$ be  of the homotopy
type of a $G-$CW complex. Let $X^k$ denote the $k-$skeleton of
$X$.  Consider the functors
$$\mbox{Map}_G(-, R_{\sigma}Y)\mbox{ and }\mbox{Map}_{C_G(\sigma)}((-)^{\sigma}, Y)$$ from
$G\mathcal{T}$ to $\mathcal{T}$. Both of them sends homotopy
colimit to  homotopy limit. In addition, we have a natural map
from $\mbox{Map}_G(-, R_{\sigma}Y)$ to $\mbox{Map}_{C_G(\sigma)}((-)^{\sigma}, Y)$ by
sending a $G-$map $F: X\longrightarrow R_{\sigma}Y$ to the composition
\begin{equation}X^{\sigma}\buildrel{F^{\sigma}}\over\longrightarrow (R_{\sigma}Y)^{\sigma}\longrightarrow
Y^{\sigma}\subseteq Y\label{fi}\end{equation} with the second map
$f\mapsto f(e)$. Note that for any $f\in (R_{\sigma}Y)^{\sigma}$, $i=1, \cdots n$,  $f(e)=(\sigma_i\cdot
f)(e)= f(e \sigma_i)=f(\sigma_i)=\sigma_i\cdot f(e)$ so $f(e)\in (Y\ast S_{C_G(\sigma),
\sigma})^{\sigma}=Y^{\sigma}$ and the second map is well-defined. It gives weak
equivalence on orbits, as shown in (\ref{weakorbitbasic}). Thus, 
$R_{\sigma}$ is a homotopical right adjoint of
$L_{\sigma}$. \end{proof}

The subgroup $\{[(1, t)]\in\Lambda_G(\sigma)|
t\in\mathbb{R}^n\}$ of $\Lambda_G(\sigma)$ is isomorphic to $\mathbb{R}^n$.
We use the same symbol $\mathbb{R}^n$ to denote it.
\begin{theorem}Let $Y$ be a $\Lambda_G(\sigma)-$space. Consider the
functor $\mathcal{L}_{\sigma}: G\mathcal{T}\longrightarrow
\Lambda_G(\sigma)\mathcal{T}, X\mapsto X^{\sigma}$ where $\Lambda_G(\sigma)$ acts
on $X^{\sigma}$ by $[g, t]\cdot x=gx.$ The functor $\mathcal{R}_{\sigma}:
\Lambda_G(\sigma)\mathcal{T}\longrightarrow G\mathcal{T}$ with
\begin{equation}\mathcal{R}_{\sigma}Y=\mbox{Map}_{C_G(\sigma)}(G, Y^{\mathbb{R}^{n}}\ast S_{C_G(\sigma), \sigma})\end{equation} is a homotopical right adjoint of $\mathcal{L}_{\sigma}$.\label{KUM}
\end{theorem}

\begin{proof}
Let $X$ be a $G-$space. Let $H$ be any closed subgroup of $G$.
For any $G-$space $X$, $\mathbb{R}^n$ acts trivially on $X^{\sigma}$,
thus, the image of any $\Lambda_G(\sigma)-$equivariant map
$X^{\sigma}\longrightarrow Y$ is in $Y^{\mathbb{R}^n}$. So we have
$\mbox{Map}_{\Lambda_G(\sigma)}(X^{\sigma}, Y)=\mbox{Map}_{C_G(\sigma)}(X^{\sigma},
Y^{\mathbb{R}^n}).$

First we show $f: (G/H)^{\sigma}\longrightarrow Y^{\mathbb{R}^n}$ extends
uniquely up to $C_G(\sigma)-$homotopy to a $C_G(\sigma)-$equivariant map
$\widetilde{f}: G/H\longrightarrow Y^{\mathbb{R}^n}\ast S_{C_G(\sigma),
\sigma}$. $f$ can be viewed as a map $(G/H)^{\sigma}\longrightarrow
Y^{\mathbb{R}^n}\ast S_{C_G(\sigma), \sigma}$ by composing with the inclusion
as the end of the join
$$Y^{\mathbb{R}^n}\longrightarrow Y^{\mathbb{R}^n}\ast S_{C_G(\sigma), \sigma},\mbox{    } y\mapsto (1y, 0).$$

The rest of the proof is analogous to that of Theorem
\ref{main}.\end{proof}

Theorem \ref{KUM} implies Theorem \ref{Emain} directly.
\begin{theorem}
For any compact Lie group $G$ and any integer $n$ and $m$, let $E_{G, n, m}$
denote the space representing the $m-$th $G-$equivariant
$E_n-$theory. Then the theory $QE^m_{n,G}$ is weakly represented  by the
space
$$QE_{G,
n, m}:=\prod_{\sigma\in G^{n}_{z}}\mathcal{R}_{\sigma}(E_{\Lambda_G(\sigma),
n, m})
$$ in the sense of  (\ref{weaksensp})
\begin{equation}\pi_0(QE_{G, n, m})=QE^m_{n, G}(S^0).\label{weaksensp}\end{equation} where
$\mathcal{R}_{\sigma}(E_{\Lambda_G(\sigma), n, m})$ is the space
$$\mbox{Map}_{C_G(\sigma)}(G, E_{\Lambda_G(\sigma), n, m}^{\mathbb{R}^n}\ast
S_{C_G(\sigma), \sigma}).$$

\label{Emain}
\end{theorem}

\section{Orthogonal $G-$spectrum of $QE^*_{n, G}$}\label{GorthospectraQEll}
In this section, we consider equivariant cohomology
theories $E^*_{G}$ that  have the same key features as equivariant
complex K-theories. More explicitly, \\ $\bullet$ The theories
$\{E_{G}^*\}_G$ have the change-of-group isomorphism, i.e. for any
closed subgroup $H$ of $G$ and $H-$space $X$, the change-of-group
map $\rho^G_H: E^*_{G}(G\times_HX)\longrightarrow E^*_{H}(X)$ defined
by $E^*_{G}(G\times_HX)\buildrel{\phi^*}\over\longrightarrow
E^*_{H}(G\times_H X)\buildrel{i^*}\over\longrightarrow E_{H}^*(X)$ is
an isomorphism where $\phi^*$ is the restriction map and $i:
X\longrightarrow
G\times_HX$ is the $H-$equivariant map defined by $i(x)=[e, x].$.\\
$\bullet$ There exists an orthogonal spectrum $E$ such that for
any compact Lie group $G$ and "large" real $G-$representation $V$
and a compact $G-$space $B$ we have a bijection
$E^{V}_{G}(B)\longrightarrow [B_+, E(V)]^G$. And $(E_{G}, \eta^{E}, \mu^{E})$ is the underlying orthogonal $G-$spectrum of $E$.\\
$\bullet$ Let $G$ be a compact Lie group and $V$ an orthogonal
$G-$representation. For every ample $G-$representation $W$, the
adjoint structure map $\widetilde{\sigma}^E_{V, W}:
E(V)\longrightarrow \mbox{Map}(S^W, E(V\oplus W))$ is a $G-$weak
equivalence.

In this section  we
construct a $\mathcal{I}_G-$FSP $(QE_n(G, -), \eta^{QE_n}, \mu^{QE_n})$ representing the theory
$QE^*_{n, G}(-)$ in the category $GwS$ defined in Definition \ref{gwsdef}.

\subsection{The construction of $QE_n(G, -)$}\label{orthonsp}

\subsubsection{The construction of $S(G, V)_{\sigma}$}

In this section, for each $\sigma\in G^n_z$, we construct an orthogonal version $S(G, V)_{\sigma}:=
Sym(V)\setminus Sym(V)^{\sigma}$ of the space $S_{G, \sigma}$. It is the space
classified by the condition (\ref{sfnew}) which is also the  condition classifying $S_{G,\sigma}$.

Let $V$ be a real $G-$representation. Let
$Sym^n(V)$ denote the $n-$th symmetric power  $V^{\otimes n} $,
which has an evident $G\wr\Sigma_n-$action on it. Let
$$Sym(V):=\bigoplus_{n\geq 0} Sym^n(V).$$

If $V$ is an ample $G-$representation, $Sym(V)$ is a faithful
$H-$representation, thus, a complete $H-$universe.

Define \begin{equation}S(G, V)_{\sigma}:=Sym(V)\setminus Sym(V)^{\sigma}=\bigcup_{i=1}^n Sym(V)\setminus Sym(V)^{\sigma_i}.\end{equation} The complex conjugation on $V$ induces an involution on it. Note that for any
subgroup $H$ of $G$ containing $\Gamma=\langle\sigma_1, \cdots \sigma_n\rangle$, $S(H, V)_\sigma$ has the same
underlying space as $S(G, V)_\sigma$.
\begin{proposition}
Let $V$ be an orthogonal $G-$representation. For any closed
subgroup $H\leqslant C_G(\sigma)$, $S(G, V)_{\sigma}$ satisfies
\begin{equation}S(G, V)_{\sigma}^H\simeq\begin{cases}\mbox{pt},&\text{if $\Gamma \nleq
H$;}\\
\emptyset,&\text{if $\Gamma\leqslant
H$.}\end{cases}\end{equation}\label{sfnew}
\end{proposition}
\begin{proof}
If $\Gamma\leqslant H$, $Sym(V)^H$ is a subspace of
$Sym(V)^{\sigma}$, so $(Sym(V)\setminus Sym(V)^{\sigma})^H$ is empty. To simplify the
symbol, we use $Sym^{n, \perp}$ denote the orthogonal complement of
$Sym^n(V)^{\sigma}$ in $Sym^n(V)$.
\[(Sym(V)\setminus Sym(V)^{\sigma})^H= colim_{n\longrightarrow\infty} (Sym^n(V)^{\sigma})^H\times
\big((Sym^{n, \perp})^H\setminus\{0\} \big)\] Then $ (Sym^{n,
\perp})^H\setminus\{0\}\backsimeq S^{k_n-1}$ where $k_n$
is the dimension of $(Sym^{n, \perp})^H$.  As $n$ goes to
infinity, $k_n$ goes to infinity. When $k_n$ is large enough,
$S^{k_n-1}$ is contractible. So $(Sym(V)\setminus Sym(V)^{\sigma})^H$ is
contractible.
\end{proof}

\subsubsection{The construction of $F_{\sigma}(G, V)$}
Next, we construct a space $F_{\sigma}(G,
V)$ representing the theory $E^{V^{\sigma}}_{\Lambda_G(\sigma)}(-)$.

If $V$ is a faithful $G-$representation,
by Proposition \ref{farithreallambda}, we have the
faithful $\Lambda_G(\sigma)-$representation $(V)^{\mathbb{R}}_{\sigma}$. In addition, $V^{\sigma}$ can be
considered as a
$\Lambda_G(\sigma)-$representation with trivial $\mathbb{R}-$action.
The space $E((V)^{\mathbb{R}}_{\sigma}\oplus V^{\sigma})$ represents
$E^{(V)^{\mathbb{R}}_{\sigma}\oplus V^{\sigma}}_{\Lambda_G(\sigma)}(-)$. So we have
$$\mbox{Map}(S^{(V)^{\mathbb{R}}_\sigma}, E((V)^{\mathbb{R}}_\sigma\oplus
V^{\sigma}))$$  represents $E^{V^\sigma}_{\Lambda_G(\sigma)}(-)$ since $$[X^\sigma,
\mbox{Map}(S^{(V)^{\mathbb{R}}_\sigma}, E((V)^{\mathbb{R}}_\sigma\oplus
V^\sigma))]^{\Lambda_G(\sigma)} $$ is isomorphic to $$[X^\sigma\wedge S^{(V)^{\mathbb{R}}_\sigma},
E((V)^{\mathbb{R}}_\sigma\oplus V^\sigma)]^{\Lambda_G(\sigma)}=
E^{(V)^{\mathbb{R}}_\sigma\oplus V^\sigma}_{\Lambda_G(\sigma)}(X^\sigma\wedge
S^{(V)^{\mathbb{R}}_\sigma})= E^{V^\sigma}_{\Lambda_G(\sigma)}(X^\sigma).$$

To simplify the symbol, we use $F_{\sigma}(G, V)$ to denote the space
$$\mbox{Map}_{\mathbb{R}}(S^{(V)^{\mathbb{R}}_\sigma},
E((V)^{\mathbb{R}}_{\sigma}\oplus V^{\sigma})).$$ Its basepoint $c_0$ is the
constant map to the basepoint of
$E((V)^{\mathbb{R}}_{\sigma}\oplus V^{\sigma})$.

$F_{\sigma}: (G, V)\mapsto F_{\sigma}(G, V)$
provides a functor from $\mathcal{I}_{G}$ to the category
$C_G(\sigma)\mathcal{T}$ of $C_G(\sigma)-$spaces. It has the properties below.
\begin{proposition}
Let $G$ and $H$ be compact Lie groups. Let $V$ be a real
$G-$representation and $W$ a real $H-$representation. Let $\sigma\in
G^{n}_z$, $\tau\in H^{n}_z$.\\ (i)  We have the unit
map $\eta_{\sigma}(G, V): S^{V^\sigma}\longrightarrow F_\sigma(G, V)$ and the
multiplication
$$\mu^F_{(\sigma, \tau)}((G, V), (H, W)): F_{\sigma}(G, V)\wedge F_{\tau}(H,
W)\longrightarrow F_{(\sigma, \tau)}(G\times H, V\oplus W)$$  making the
unit, associativity and centrality of unit diagram  commute. And
$\eta_\sigma(G, V)$ is $C_G(\sigma)-$equivariant and $\mu^F_{(\sigma, \tau)}((G, V),
(H, W))$ is $C_{G\times H}(\sigma, \tau)-$ equivariant. \\ (ii)Let
$\Delta_G$ denote the diagonal map $G\longrightarrow G\times
G,\mbox{      }g\mapsto (g, g).$ Let $\widetilde{\sigma}_{\sigma}(G, V,
W): F_{\sigma}(G, V)\longrightarrow \mbox{Map}(S^{W^{\sigma}}, F_{\sigma}(G, V\oplus
W))$ denote the map
$$x\mapsto (w\mapsto \big(\Delta^*_{G}\circ \mu^F_{(\sigma, \tau)}((G, V), (G, W))\big)\big(x, \eta_{\sigma}(G, W)(w)\big)).$$ Then $\widetilde{\sigma}_\sigma(G, V, W)$
is a $\Lambda_G(\sigma)-$weak equivalence when $V$ is an ample
$G-$representation. \\ (iii) If $(E, \eta^{E}, \mu^{E})$ is
commutative, we have
\begin{equation}\mu^F_{(\sigma, \tau)}((G, V), (H, W))(x\wedge y)=
\mu^F_{(\tau, \sigma)}((H, W), (G, V))(y\wedge x)
\label{quasicommFg}\end{equation} for any $x\in F_{\sigma}(G, V)$ and
$y\in F_{\tau}(H, W)$.
\label{newF}
\end{proposition} 

The proof is straightforward and left to the readers.

\subsubsection{The construction of $QE_n(G, V)$}
\label{qegcon}Recall in Theorem \ref{Emain} we  construct a $G-$space
$QE_{G, n, m}$ representing $QE^m_{n, G}(-)$. In this section we go a step further.

Apply Theorem \ref{KUM}, we get the conclusion below.
\begin{proposition}
Let $V$ be a faithful orthogonal $G-$representation. Let $B'_n(G,
V)$ denote the space $$\prod_{\sigma\in
G^{n}_{z}}\mbox{Map}_{C_G(\sigma)}(G, F_\sigma(G, V)\ast S(G, V)_\sigma).$$
$QE^V_{n, G}(-)$ is weakly represented by $B'_n(G, V)$
in the sense $\pi_0(B'_n(G, V))= QE^V_{n, G}(S^0). $\label{quickEFS}
\end{proposition}
The proof of Proposition \ref{quickEFS} is analogous to that of
Theorem \ref{Emain} step by step.

Below is the main theorem in Section \ref{orthonsp}. We will
use formal linear combination $$t_1a+t_2b\mbox{ with }0\leqslant
t_1, t_2\leqslant 1, t_1+t_2=1$$ to denote points in join.

\begin{proposition}Let $QE_{n, \sigma}(G, V)$ denote $$\{t_1a+t_2b\in F_\sigma(G, V)\ast S(G, V)_{\sigma}| \|b\|\leqslant t_2\}/\{t_1c_0+t_2
b\}.
$$ It is the quotient space of a closed subspace of  the join $F_{\sigma}(G, V)\ast S(G,
V)_{\sigma}$ with all the points of the form $t_1c_0+t_2 b$ collapsed to
one point, which we pick as the basepoint of $QE_{n, \sigma}(G, V)$, where
$c_0$ is the basepoint of $F_\sigma(G, V)$. $QE_{n, \sigma}(G, V)$ has the
evident $C_G(\sigma)-$action. And it is $C_G(\sigma)-$weak equivalent to
$F_\sigma(G, V)\ast S(G, V)_\sigma$. As a result, $\prod\limits_{\sigma\in
G^{n}_{z}}\mbox{Map}_{C_G(\sigma)}(G, QE_{n, \sigma}(G, V))$ is $G-$weak
equivalent to $\prod\limits_{\sigma\in
G^{n}_{z}}\mbox{Map}_{C_G(\sigma)}(G, F_\sigma(G, V)\ast S(G, V)_\sigma)$.
So when $V$ is a faithful $G-$representation,
\begin{equation}QE_n(G, V):=\prod\limits_{\sigma\in
G^{n}_{z}}\mbox{Map}_{C_G(\sigma)}(G, QE_{n, \sigma}(G, V))
\label{EGVnog}\end{equation} weakly represents $QE^V_{n, G}(-)$ in the
sense $\pi_0(QE_n(G, V))\cong QE^V_{n, G}(S^0).$\label{anoQEllwere}
\end{proposition}
\begin{proof}
First we show $F_\sigma(G, V)\ast S(G, V)_\sigma$ is $C_G(\sigma)-$homotopy
equivalent to $$QE'_{n, \sigma}(G, V):=\{t_1a+t_2b\in F_\sigma(G, V)\ast S(G,
V)_\sigma| \|b\|\leqslant t_2\}.$$

Note that $b\in S(G, V)_\sigma$ is never zero. Let $j: QE'_{n, \sigma}(G,
V)\longrightarrow F_\sigma(G, V)\ast S(G, V)_{\sigma}$ be the inclusion. Let
$p: F_{\sigma}(G, V)\ast S(G, V)_{\sigma}\longrightarrow QE'_{n, \sigma}(G, V)$ be the
$C_G(\sigma)-$map sending $t_1a+t_2b$ to $t_1a+t_2\frac{min\{\|b\|,
t_2\}}{\|b\|}b$. Both $j$ and $p$ are both continuous and
$C_G(\sigma)-$equivariant. $p\circ j$ is the identity map of $QE'_{n, \sigma}(G,
V)$. We can define a $C_G(\sigma)-$homotopy
$$H: (F_\sigma(G, V)\ast S(G, V)_\sigma)\times I\longrightarrow F_\sigma(G,
V)\ast S(G, V)_\sigma$$ from the identity map on $F_{\sigma}(G, V)\ast S(G,
V)_{\sigma}$ to $j\circ p$ by shrinking. For any $t_1a+t_2b\in F_\sigma(G,
V)\ast S(G, V)_\sigma$, Define
\begin{equation}H(t_1a+t_2b, t):= t_1a+t_2((1-t)b+t\frac{min\{\|b\|, t_2\}}{\|b\|}b). \end{equation}

Then we show $QE'_{n, \sigma}(G, V)$ is $G-$weak equivalent to $QE_{n, \sigma}(G, V)$.
Let $q: QE'_{n, \sigma}(G, V)\longrightarrow QE_{n, \sigma}(G, V)$ be the quotient
map. Let $H$ be a closed subgroup of $C_G(\sigma)$.

If the group $\Gamma$ is in $H$, since $S(G, V)_{\sigma}^H$ is empty, so $QE_{n, \sigma}(G, V)^H$
is in the end $F_{\sigma}(G, V)$ and can be identified with $F_{\sigma}(G,
V)^H$. In this case
$q^H$ is the identity map.

If $\Gamma$ is not in $H$, $QE'_{n, \sigma}(G, V)^H$ is contractible. The cone
$\{c_0\}\ast S(G, V)_{\sigma}^H$ is contractible, so $q\big((\{c_0\}\ast
S(G, V)_{\sigma})^H\big)=q(\{c_0\}\ast S(G, V)_{\sigma}^H)$ is contractible.
Note that the subspace of all the points of the form $t_1c_0+t_2b$
for any $t_1$
and $b$  is $q\big((\{c_0\}\ast S(G, V)_{\sigma})^H\big)$. Therefore, $QE_{n, \sigma}(G, V)^H=QE'_{n, \sigma}(G, V)^H/q(\{c_0\}\ast S(G,
V)_\sigma)^H$ is contractible.

Therefore, $QE'_{\sigma}(G, V)$ is $G-$weak equivalent to $F_{\sigma}(G, V)\ast
S(G, V)_{\sigma}$.
\end{proof}

\begin{proposition}
Let $\sigma\in G^{n}_z$. Let $Y$ be a based $\Lambda_G(\sigma)-$space.  Let
$\widetilde{Y}_{\sigma}$ denote the $C_G(\sigma)-$space
$$\{t_1a+t_2b\in Y^{\mathbb{R}^n}\ast S(G, V)_{\sigma}| \|b\|\leqslant
t_2\}/\{t_1y_0+t_2 b\}.
$$ It is the quotient space of a closed subspace of $Y^{\mathbb{R}^n}\ast S(G,
V)_{\sigma}$ with all the points of the form     $t_1y_0+t_2 b$
  collapsed to one point, i.e the basepoint of $\widetilde{Y}_{\sigma}$,
where $y_0$ is the basepoint of $Y$. $\widetilde{Y}_\sigma$ is
$C_G(\sigma)-$weak equivalent to $Y^{\mathbb{R}^n}\ast S(G, V)_\sigma$. As a
result, the functor $R_\sigma: C_G(\sigma)\mathcal{T}\longrightarrow
G\mathcal{T}$ with $R_{\sigma}\widetilde{Y}=\mbox{Map}_{C_G(\sigma)}(G,
\widetilde{Y}_{\sigma})$ is a homotopical right adjoint of $L:
G\mathcal{T}\longrightarrow C_G(\sigma)\mathcal{T}\mbox{, }X\mapsto
X^{\sigma}$.

\label{generalYS}
\end{proposition}

The proof 
is analogous to that of
Theorem \ref{KUM} and Proposition \ref{anoQEllwere}.

\begin{remark}
We can consider $QE_{n, \sigma}(G, V)$  as a quotient space of a subspace of
$F_{\sigma}(G, V)\times Sym(V)\times I$  \begin{equation}\{(a, b, t)\in
F_{\sigma}(G, V)\times Sym(V)\times I|  \|b\|\leqslant t; \mbox{   and }
b\in S(G, V)_{\sigma}\mbox{   if   }t\neq
0\}\label{qimiaoanoano}\end{equation} by identifying points $(a,
b, 1)$ with $(a', b, 1)$, and collapsing all the points $(c_0, b,
t)$ for any $b$ and $t$. In other words, the end $F_{\sigma}(G, V)$ in
the join $F_{\sigma}(G, V)\ast S(G, V)_{\sigma}$ is identified with the points
of the form $(a, 0, 0)$ in (\ref{qimiaoanoano}).
\label{ananotherEll}
\end{remark}

\begin{proposition}For each $\sigma\in G^{n}_z$, $$QE_{n, \sigma}: \mathcal{I}_G\longrightarrow C_G(\sigma)\mathcal{T},\mbox{   } (G, V)\mapsto QE_{n, \sigma}(G, V)$$ is a well-defined functor.
As a result, $$QE_n: \mathcal{I}_G\longrightarrow
G\mathcal{T},\mbox{ } (G, V)\mapsto \prod\limits_{\sigma\in
G^{n}_{z}}\mbox{Map}_{C_G(\sigma)}(G, QE_{n, \sigma}(G, V))$$ is a
well-defined functor. \end{proposition}
\begin{proof}
Let $V$ and $W$ be $G-$representations and $f: V\longrightarrow
W$ a linear isometric isomorphism. Then $f$ induces a
$C_G(\sigma)-$homeomorphism $F_\sigma(f)$ from $F_\sigma(G, V)$ to $F_\sigma(G, W)$
and a $C_G(\sigma)-$homeomorphism $S_\sigma(f)$ from $S(G,V)_\sigma$  to
$S(G,W)_{\sigma}$. We have the well-defined map
$$QE_{n, \sigma}(f): QE_{n, \sigma}(G, V)\longrightarrow QE_{n, \sigma}(G, W)$$ sending a point
represented by $t_1a+t_2b$ in the join to that represented by
$t_1F_\sigma(f)(a)+t_2S_\sigma(f)(b)$. And $QE_n(f): QE_n(G, V)\longrightarrow
QE_n(G, W)$ is defined by
$$\prod\limits_{\sigma\in G^{n}_{z}} \alpha_\sigma\mapsto\prod\limits_{\sigma\in G^{n}_{z}} QE_{n, \sigma}(f)\circ\alpha_\sigma.$$

It is straightforward to check that all the axioms hold.
\end{proof}

\subsection{Construction of $\eta^{QE_n}$ and $\mu^{QE_n}$}\label{Gorthospstructure}
In this section we construct a unit map $\eta^{QE_n}$ and a
multiplication $\mu^{QE_n}$ so that we get a commutative
$\mathcal{I}_G-$FSP representing the $QE_n-$theory in $GwS$.

Let $G$ and $H$ be compact Lie groups, $V$ an orthogonal
$G-$representation and $W$ an orthogonal $H-$representation. Let $\sigma\in G^{n}_z$. We
use $x_\sigma$ to denote the basepoint of $QE_{n, \sigma}(G, V)$, which is
defined in Proposition \ref{anoQEllwere}.  For
each $v\in S^V$, there are $v_1\in S^{V^\sigma}$ and $v_2\in
S^{(V^\sigma)^{\perp}}$ such that $v=v_1\wedge v_2$. Let
$\eta^{QE_n}_{\sigma}(G, V): S^V\longrightarrow QE_{n, \sigma}(G, V)$ be the map
\begin{equation}\eta^{QE_n}_\sigma(G,
V)(v):=\begin{cases}(1-\|v_2\|)\eta_\sigma(G, V)(v_1)+ \|v_2\|v_2,
&\text{if $\|v_2\|\leqslant 1$;}\\ x_\sigma, &\text{if
$\|v_2\|\geqslant 1$.}
\end{cases}\label{finaleta}\end{equation} 
\begin{lemma}The map $\eta^{QE_n}_\sigma(G, V)$ defined in (\ref{finaleta}) is well-defined, continuous and $C_G(\sigma)-$equivariant.
\label{etaEcontinuous}\end{lemma}

\begin{remark}
For any $\sigma\in G^{n}_z$, it's straightforward to check the diagram
below commutes.
$$\begin{CD}S^{V^{\sigma}} @>\eta_{\sigma}(G, V)>> F_{\sigma}(G, V) \\
@VVV      @VVV \\
S^V @>\eta^{QE_n}_{\sigma}(G, V)>> QE_{n, \sigma}(G, V)\end{CD}$$ where both vertical
maps are inclusions. By Lemma \ref{etaEcontinuous}, the map
\begin{equation}\eta^{QE_n}(G, V): S^V\longrightarrow
\prod\limits_{\sigma\in G^{n}_{z}}\mbox{Map}_{C_G(\sigma)}(G, QE_{n, \sigma}(G,
V))\mbox{,   } v \mapsto \prod\limits_{\sigma\in
G^{n}_{z}}(\alpha\mapsto \eta^{QE_n}_{\sigma}(G, V)(\alpha\cdot
v)),\label{GetaQEll}
\end{equation}
is well-defined and continuous. Moreover, $\eta^{QE_n}:
S\longrightarrow QE_n $ with $QE_n(G, V)$ defined in (\ref{EGVnog}) is
well-defined. \end{remark} Next, we construct the multiplication
map $\mu^{QE_n}$. First we define a map $$\mu^{QE_n}_{(\sigma, \tau)}((G,V),
(H, W)): QE_{n, \sigma}(G, V)\wedge QE_{n, \tau}(H, W)\longrightarrow QE_{n, (\sigma,
\tau)}(G\times H, V\oplus W)$$ by sending a point $[t_1a_1+
t_2b_1]\wedge [u_1a_2 + u_2b_2]$ to
\begin{equation}
\begin{cases}[(1-\sqrt{t^2_2+u^2_2})\mu^F_{(\sigma, \tau)}((G,V)
, (H, W))(a_1\wedge a_2) &\text{if $t^2_2+u^2_2\leq 1$ and
$t_2u_2\neq 0$;}\\ + \sqrt{t^2_2+u^2_2}(b_1+b_2)],
& \\
[(1-t_2)\mu^F_{(\sigma, \tau)}((G,V), (H, W))(a_1\wedge a_2)+ t_2b_1],
&\text{if $u_2=0$ and $0<t_2<1$}; \\ [(1-u_2)\mu^F_{(\sigma, \tau)}((G,V), (H, W))(a_1\wedge a_2) + u_2b_2],
&\text{if $t_2=0$ and $0<u_2<1$; } \\ [1\mu^F_{(\sigma, \tau)}((G,V), (H, W))(a_1\wedge a_2) + 0],
&\text{if $u_2=0$ and $t_2=0$; }\\ x_{\sigma, \tau},
&\text{Otherwise.}\end{cases}\label{muEg}
\end{equation} where $x_{\sigma, \tau}$ is the basepoint of
$QE_{n, (\sigma, \tau)}(G\times H, V\oplus W)$.

\begin{lemma}The map $\mu^{QE_n}_{(\sigma, \tau)}((G,V), (H, W))$ defined in (\ref{muEg}) is well-defined and
continuous. \label{muEcontinuous}
\end{lemma}

The basepoint of $QE_n(G, V)$ is the product of the basepoint of
each factor $\mbox{Map}_{C_G(\sigma)}(G, QE_{n, \sigma}(G, V))$, i.e. the product
of the constant map to the basepoint of  each $QE_{n, \sigma}(G, V)$.

We can define the multiplication $\mu^{QE_n}((G, V), (H, W)): QE_n(G,
V)\wedge QE_n(H, W)\longrightarrow QE_n(G\times H, V\oplus W)$ by
$$\big(\!\!\prod\limits_{\sigma\in G^{n}_{z}}
\!\!\alpha_\sigma\big)\wedge\big(\!\!\prod\limits_{ \tau\in
H^{n}_{z}}\!\!\beta_{\tau}\big)\mapsto\prod\limits_{\substack{\sigma\in
G^{n}_{z}\\ \tau\in H^{n}_{z}}}\!\!\bigg(\!\!(\sigma',
\tau')\mapsto \mu^{QE_n}_{(\sigma, \tau)}((G, V), (H,
W))\big(\alpha_\sigma(\sigma')\wedge\beta_\tau(\tau')\big)\!\!\bigg).
$$

\begin{theorem}

$QE_n:\mathcal{I}_G\longrightarrow G\mathcal{T}$ together with the
unit map $\eta^{QE_n}$ defined in (\ref{GetaQEll}) and  the
multiplication $\mu^{QE_n}((G, -), (G, -))$ gives a
commutative $\mathcal{I}_G-$FSP that weakly represents
$QE^*_{n, G}(-)$. \label{GorthoQEll}\end{theorem}

\begin{remark}
We apply a conclusion from Chapter 3, Section 1, in \cite{SS}.  A
$G-$spectrum $Y$ is isomorphic to an orthogonal $G-$spectrum of
the form
$X\langle G\rangle$ for some orthogonal spectrum $X$ if and only if for every trivial $G-$representation $V$ the $G-$action on
$Y(V)$ is trivial. $QE_n(G, V)$ is not trivial when $V$ is trivial. So
it cannot arise from an orthogonal spectrum. \label{Enotorthosp}
\end{remark}

\begin{proposition}
Let $G$ be any compact Lie group. Let $V$ be an ample orthogonal
$G-$representation and $W$ an orthogonal $G-$representation. Let
$\sigma^{QE_n}_{G, V, W}: S^W\wedge QE_n(G, V)\longrightarrow QE_n(G,
V\oplus W)$ denote the structure map of $QE_n$ defined by the unit
map  $\eta^{QE_n}(G, V)$. Let $\widetilde{\sigma}^{QE_n}_{G, V, W}$
denote the right adjoint of $\sigma^{QE_n}_{G, V, W}$. Then
$\widetilde{\sigma}^{QE_n}_{G, V, W}: QE_n(G, V)\longrightarrow
\mbox{Map}(S^W, QE_n(G, V\oplus W))$ is a $G-$weak equivalence.
\label{fibrantweakequiv}\end{proposition}

Let $G$ and $H$ be compact Lie groups, $V$ an orthogonal
$G-$representation and $W$ an orthogonal $H-$representation. We
use $x_\sigma$ to denote the basepoint of $QE_\sigma(G, V)$, which is
defined in Proposition \ref{anoQEllwere}. Let $\sigma\in G^{n}_z$. For
each $v\in S^V$, there are $v_1\in S^{V^\sigma}$ and $v_2\in
S^{(V^\sigma)^{\perp}}$ such that $v=v_1\wedge v_2$. Let
$\eta^{QE_n}_{\sigma}(G, V): S^V\longrightarrow QE_{n, \sigma}(G, V)$ be the map
\begin{equation}\eta^{QE_n}_{\sigma}(G,
V)(v):=\begin{cases}(1-\|v_2\|)\eta_\sigma(G, V)(v_1)+ \|v_2\|v_2,
&\text{if $\|v_2\|\leqslant 1$;}\\ x_\sigma, &\text{if
$\|v_2\|\geqslant 1$.}
\end{cases}\label{finaleta}\end{equation}
The map $\eta^{QE_n}_{\sigma}(G, V)$ defined in (\ref{finaleta}) is well-defined, continuous and $C_G(\sigma)-$equivariant.

\begin{remark}
For any $\sigma\in G^{n}_z$, it's straightforward to check the diagram
below commutes.
$$\begin{CD}S^{V^\sigma} @>\eta_{\sigma}(G, V)>> F_{\sigma}(G, V) \\
@VVV      @VVV \\
S^V @>\eta^{QE_n}_{\sigma}(G, V)>> QE_{n, \sigma}(G, V)\end{CD}$$ where both vertical
maps are inclusions. By Lemma \ref{etaEcontinuous}, the map
\begin{equation}\eta^{QE_n}(G, V): S^V\longrightarrow
\prod\limits_{\sigma\in G^{n}_{z}}\mbox{Map}_{C_G(\sigma)}(G, QE_{n, \sigma}(G,
V))\mbox{,   } v \mapsto \prod\limits_{\sigma\in
G^{n}_{z}}(\alpha\mapsto \eta^{QE_n}_\sigma(G, V)(\alpha\cdot
v)),\label{GetaQEll}
\end{equation}
is well-defined and continuous. Moreover, $\eta^{QE_n}:
S\longrightarrow QE_n$ with $QE_n(G, V)$ defined in (\ref{EGVnog}) is
well-defined. \end{remark} Next, we construct the multiplication
map $\mu^{QE_n}$. First we define a map $$\mu^{QE_n}_{(\sigma, \tau)}((G,V),
(H, W)): QE_{n, \sigma}(G, V)\wedge QE_{n, \tau}(H, W)\longrightarrow QE_{n, (\sigma, \tau)}(G\times H, V\oplus W)$$ by sending a point $[t_1a_1+
t_2b_1]\wedge [u_1a_2 + u_2b_2]$ to
\begin{equation}
\begin{cases}[(1-\sqrt{t^2_2+u^2_2})\mu^F_{(\sigma, \tau)}((G,V)
, (H, W))(a_1\wedge a_2) &\text{if $t^2_2+u^2_2\leq 1$ and
$t_2u_2\neq 0$;}\\ + \sqrt{t^2_2+u^2_2}(b_1+b_2)],
& \\
[(1-t_2)\mu^F_{(\sigma, \tau)}((G,V), (H, W))(a_1\wedge a_2)+ t_2b_1],
&\text{if $u_2=0$ and $0<t_2<1$}; \\ [(1-u_2)\mu^F_{(\sigma, \tau)}((G,V), (H, W))(a_1\wedge a_2) + u_2b_2],
&\text{if $t_2=0$ and $0<u_2<1$; } \\ [1\mu^F_{(\sigma, \tau)}((G,V), (H, W))(a_1\wedge a_2) + 0],
&\text{if $u_2=0$ and $t_2=0$; }\\ x_{\sigma, \tau},
&\text{Otherwise.}\end{cases}\label{muEg}
\end{equation} where $x_{\sigma, \tau}$ is the basepoint of
$QE_{n, (\sigma, \tau)}(G\times H, V\oplus W)$.
The map $\mu^{QE_n}_{(\sigma, \tau)}((G,V), (H, W))$ defined in (\ref{muEg}) is well-defined and
continuous.

The basepoint of $QE_n(G, V)$ is the product of the basepoint of
each factor $\mbox{Map}_{C_G(\sigma)}(G, QE_{n, \sigma}(G, V))$, i.e. the product
of the constant map to the basepoint of  each $QE_{n, \sigma}(G, V)$.

We can define the multiplication $\mu^{QE_n}((G, V), (H, W)): QE_n(G,
V)\wedge QE_n(H, W)\longrightarrow QE_n(G\times H, V\oplus W)$ by
$$\big(\!\!\prod\limits_{\sigma\in G^{n}_{z}}
\!\!\alpha_\sigma\big)\wedge\big(\!\!\prod\limits_{ \tau\in
H^{n}_{z}}\!\!\beta_\tau\big)\mapsto\prod\limits_{\substack{\sigma\in
G^{n}_{z}\\ \tau\in H^{n}_{z}}}\!\!\bigg(\!\!(\sigma',
\tau')\mapsto \mu^{QE_n}_{(\sigma, \tau)}((G, V), (H,
W))\big(\alpha_\sigma(\sigma')\wedge\beta_\tau(\tau')\big)\!\!\bigg).
$$

\begin{theorem}

$QE_n(G, -):\mathcal{I}_G\longrightarrow G\mathcal{T}$ together with the
unit map $\eta^{QE_n}$ defined in (\ref{GetaQEll}) and  the
multiplication $\mu^{QE_n}((G, -), (G, -))$ gives a
commutative $\mathcal{I}_G-$FSP that weakly represents
$QE^*_{n, G}(-)$. \label{GorthoQEll}\end{theorem}
The proof of Theorem \ref{GorthoQEll} is analogous to that of Theorem 6.12 \cite{HuanSpec}.
\begin{remark}
We apply a conclusion from Chapter 3, Section 1, in \cite{SS}.  A
$G-$spectrum $Y$ is isomorphic to an orthogonal $G-$spectrum of
the form
$X\langle G\rangle$ for some orthogonal spectrum $X$ if and only if for every trivial $G-$representation $V$ the $G-$action on
$Y(V)$ is trivial. $QE_n(V)$ is not trivial when $V$ is trivial. So
it cannot arise from an orthogonal spectrum. \label{Enotorthosp}
\end{remark}

In addition, we have the conclusion below.
\begin{proposition}
Let $G$ be any compact Lie group. Let $V$ be an ample orthogonal
$G-$representation and $W$ an orthogonal $G-$representation. Let
$\sigma^{QE_n}_{G, V, W}: S^W\wedge QE_n(G, V)\longrightarrow QE_n(G,
V\oplus W)$ denote the structure map of $QE_n$ defined by the unit
map  $\eta^{QE_n}(G, V)$. Let $\widetilde{\sigma}^{QE_n}_{G, V, W}$
denote the right adjoint of $\sigma^{QE_n}_{G, V, W}$. Then
$\widetilde{\sigma}^{QE_n}_{G, V, W}: QE_n(G, V)\longrightarrow
\mbox{Map}(S^W, QE_n(G, V\oplus W))$ is a $G-$weak equivalence.
\label{fibrantweakequiv}\end{proposition}

The proof is analogous to that of Proposition 6.14 \cite{HuanSpec}.

At last, we get the main conclusion of Section
\ref{GorthospectraQEll}.
\begin{theorem}
For each positive integer $n$ and each compact Lie group $G$, there is a well-defined functor $\mathcal{Q}_{G, n}$ from the category of orthogonal ring spectra to the category of
$\mathcal{I}_G-$FSP sending $E$  to
$(QE_n(G, -), \eta^{QE_n}, \mu^{QE_n})$ that weakly represents the
quasi-theory $QE_{n, G}^*$.

\end{theorem}

\appendix

\section{Faithful representation of $\Lambda_G(\sigma)$}\label{reallambda}
We
discuss complex and real $\Lambda_G(\sigma)-$representations in
Section \ref{prere} and \ref{realdef} respectively.

\subsection{Preliminaries: faithful representations of
$\Lambda_{G}(\sigma)$}\label{prere}  In this
section, we construct a faithful
$\Lambda_G(\sigma)-$representation from a faithful
$G-$representation.

Let $G$ be a compact Lie group and  $\sigma\in G^{n}_z$. Let $l_i$ denote
the order of $\sigma_i$. Let $\rho$ denote a complex $G-$representation with
underlying space $V$. Let $i: C_G(\sigma)\hookrightarrow G$ denote
the inclusion.
Let $\{\lambda\}$ denote all the irreducible complex
representations of $C_G(\sigma)$. As said in \cite{FH}, we have
the decomposition of a representation into its isotypic components
$i^*V\cong\bigoplus\limits_{\lambda}V_{\lambda}$ where $V_{\lambda}$
denotes the sum of all subspaces of $V$ isomorphic to $\lambda$.
Each $V_{\lambda}=Hom_{C_G(\sigma)}(\lambda,
V)\otimes_{\mathbb{C}}\lambda$ is unique as a subspace. Note that each $\sigma_i$ acts on each $V_{\lambda}$ as a diagonal matrix.

Each $V_{\lambda}$ can be equipped with a
$\Lambda_G(\sigma)-$action. Each $\lambda(\sigma_i)$ is of the form
$e^{\frac{2\pi i m_{\lambda i} }{l_i}}I$ with $0< m_{\lambda i}\leq l_i$
and $I$ the identity matrix. As shown in Lemma \ref{cl},
we have the well-defined complex
$\Lambda_G(\sigma)-$representations
\[(V_{\lambda})_{\sigma}:=V_{\lambda}\odot_{\mathbb{C}}
(q^{\frac{m_{\lambda 1}}{l_1}}\otimes\cdots\otimes q^{\frac{m_{\lambda n}}{l_n}})\] and
\begin{equation}(V)_{\sigma}:=\bigoplus_{\lambda}(V_{\lambda})_{\sigma}.\label{vgrepresentation}\end{equation}

\begin{proposition}
Let $V$ be a faithful $G-$representation. Let $\sigma\in
G^{n}_z$.

(i) $(V)_{\sigma}\oplus (V)_{\sigma}\otimes_{\mathbb{C}}q^{-1}$
is a faithful $\Lambda_G(\sigma)-$representation.

(ii) $(V)_{\sigma}\oplus V^{\sigma}$ is a faithful
$\Lambda_G(\sigma)-$representation.

\label{farithlambda}
\end{proposition}

\begin{proof}

(i) Let $[a, t]\in \Lambda_G(\sigma)$ be an element acting
trivially on $V_{\sigma}$. Consider the subrepresentations
$(V_{\lambda})_{\sigma}$ and
$(V_{\lambda})_{\sigma}\otimes_{\mathbb{C}}q^{-1}$ of
$(V)_{\sigma}\oplus (V)_{\sigma}\otimes_{\mathbb{C}}q^{-1}$
respectively.  Let $v$ be an element in the underlying vector
space  $V_{\lambda}$. On $(V_{\lambda})_{\sigma}$, $[a, t]\cdot
v=e^{2\pi i t(\frac{ m_{\lambda 1}}{l_1}+\cdots \frac{ m_{\lambda n}}{l_n})}a\cdot v= v$; and on
$(V_{\lambda})_{\sigma}\otimes_{\mathbb{C}}q^{-1}$, $[a, t]\cdot
v=e^{2\pi i t(\frac{ m_{\lambda 1}}{l_1}+\cdots \frac{ m_{\lambda n}}{l_n})-2\pi i t}a\cdot v=v$. So we get
$e^{2\pi it}\cdot v=v$. Thus, $t=0$. $C_G(\sigma)$ acts faithfully
on $V$, so it acts faithfully on $(V)_{\sigma}\oplus
(V)_{\sigma}\otimes_{\mathbb{C}}q^{-1}$. Since $[a, 0]\cdot w=w$,
for any $w\in (V)_{\sigma}\oplus
(V)_{\sigma}\otimes_{\mathbb{C}}q^{-1}$, so $a=e$.

Thus, $(V)_{\sigma}\oplus (V)_{\sigma}\otimes_{\mathbb{C}}q^{-1}$
is a faithful $\Lambda_G(\sigma)-$representation.

(ii)  Note that $V^{\sigma}$ with the trivial $\mathbb{R}-$action
is the representation
$(V^{\sigma})_{\sigma}\otimes_{\mathbb{C}}q^{-1}$. The
representation $(V)_{\sigma}\oplus V^{\sigma}$ contains a
subrepresentation $(V^{\sigma})_{\sigma}\oplus
(V^{\sigma})_{\sigma}\otimes_{\mathbb{C}}q^{-1}$, which is a
faithful $\Lambda_G(\sigma)-$representation by Proposition \ref{farithlambda} (i). So
$(V)_{\sigma}\oplus V^{\sigma}$ is faithful.\end{proof}

\begin{lemma}
For any $\sigma\in G^{n}_z$,  $(-)_{\sigma}$ defined in
(\ref{vgrepresentation}) is a functor from the category of
$G-$spaces to the category of $\Lambda_G(\sigma)-$spaces.
Moreover, $(-)_{\sigma}\oplus (-)_{\sigma}\otimes_{\mathbb{C}}
q^{-1}$ and $(-)_{\sigma}\oplus (-)^{\sigma}$ in Proposition
\ref{farithlambda} are also well-defined functors from the
category of $G-$spaces to the category of
$\Lambda_G(\sigma)-$spaces.
\end{lemma}

\begin{proof}
Let $f: V\longrightarrow W$ be a $G-$equivariant map. Then $f$ is
$C_G(\sigma)-$equivairant for each $\sigma\in G^{n}_z$. For each
irreducible complex $C_G(\sigma)-$representation $\lambda$, $f:
V_{\lambda}\longrightarrow W_{\lambda}$ is
$C_G(\sigma)-$equivairant. And $f_{\sigma}:
(V_{\lambda})_{\sigma}\longrightarrow
(W_{\lambda})_{\sigma},\mbox{  } v\mapsto f(v)$ with the same
underlying spaces is well-defined and is
$\Lambda_G(\sigma)-$equivariant. It is straightforward to check if
we have two $G-$equivariant maps $f: V\longrightarrow W$ and $g:
U\longrightarrow V$, then $(f\circ g)_{\sigma}=f_{\sigma}\circ
g_{\sigma}.$ So $(-)_{\sigma}$ gives a well-defined functor from
the category of $G-$representations to the category of
$\Lambda_G(\sigma)-$representation.

The other conclusions can be proved in a similar way.\end{proof}
\begin{proposition} Let $H$ and $G$ be two compact Lie groups. Let $\sigma \in G^n_z$ and $\tau\in H^n_z$. Let $V$ be a $G-$representation and $W$ a $H-$representation.

(i) We have the isomorphisms of representations $(V\oplus
W)_{(\sigma, \tau)}=(V_{\sigma}\oplus
                     W_{\tau})$ as $\Lambda_{G\times H}(\sigma, \tau)\cong
                     \Lambda_G(\sigma)\times_{\mathbb{T}^n}\Lambda_H(\tau)-$representations;\\
$(V\oplus W)_{(\sigma, \tau)}\oplus (V\oplus W)_{(\sigma,
\tau)}\otimes_{\mathbb{C}}q^{-1}=((V)_{\sigma}\oplus
                    (V)_{\sigma}\otimes_{\mathbb{C}}q^{-1})\oplus ((W)_{\tau}\oplus
                     (W)_{\tau}\otimes_{\mathbb{C}}q^{-1})$ as $\Lambda_{G\times H}(\sigma, \tau)\cong
                     \Lambda_G(\sigma)\times_{\mathbb{T}^n}\Lambda_H(\tau)-$representations;\\
and  $(V\oplus W)_{(\sigma, \tau)}\oplus (V\oplus W)^{(\sigma,
\tau)}=((V)_{\sigma}\oplus V^{\sigma})\oplus ((W)_{\tau}\oplus
W^{\tau})$ as $\Lambda_{G\times H}(\sigma, \tau)\cong
                     \Lambda_G(\sigma)\times_{\mathbb{T}^n}\Lambda_H(\tau)-$representations.

(ii) Let $\phi: H\longrightarrow G$ be a group homomorphism. Let
$\phi_{\tau}: \Lambda_H(\tau)\longrightarrow
\Lambda_G(\phi(\tau))$ denote the group homomorphism obtained from
$\phi$. Then we have
$$\phi_{\tau}^*(V)_{\phi(\tau)}=(V)_{\tau},$$
$$\phi_{\tau}^*((V)_{\phi(\tau)}\oplus
(V)_{\phi(\tau)}\otimes_{\mathbb{C}}q^{-1})=(V)_{\tau}\oplus
(V)_{\tau}\otimes_{\mathbb{C}}q^{-1},$$
$$\phi_{\tau}^*((V)_{\phi(\tau)}\oplus V^{\phi(\tau)})=(V)_{\tau}\oplus V^{\tau}$$ as
                    $\Lambda_H(\tau)-$representations.
\end{proposition}
\begin{proof}
(i) Let $\{\lambda_G\}\mbox{    and    }\{\lambda_H\}$ denote the
sets of all the irreducible $C_G(\sigma)-$\\ representations and
all the irreducible $C_H(\tau)-$representations. Then $\lambda_G$
and $\lambda_H$ are irreducible representations of $C_{G\times
H}(\sigma, \tau)$ via the inclusion $C_G(\sigma)\longrightarrow
C_{G\times H}(\sigma, \tau)$ and $C_H(\tau)\longrightarrow
C_{G\times H}(\sigma, \tau)$.

The $\mathbb{R}-$representation assigned to each $C_{G\times
H}(\sigma, \tau)-$irreducible representation in $V\oplus W$ is the
same as that assigned to the irreducible representations of $V$
and $W$. So we have $$(V\oplus W)_{(\sigma,
\tau)}=(V_{\sigma}\oplus
                     W_{\tau})$$ as $\Lambda_{G\times H}(\sigma, \tau)\cong
                     \Lambda_G(\sigma)\times_{\mathbb{T}^n}\Lambda_H(\tau)-$representations.

Similarly we can prove the other two conclusions in (i).

(ii) Let $\sigma=\phi(\tau)$. If $(\phi_{\tau}^*V)_{\lambda_H}$ is
a $C_H(\tau)-$subrepresentation of $\phi_{\tau}^*V_{\lambda_G}$,
the $\mathbb{R}-$representation assigned to it is the same as that
to $V_{\lambda_G}$. So we have
$\phi_{\tau}^*(V)_{\phi(\tau)}=(V)_{\tau}$ as
$\Lambda_H(\tau)-$representations.

Similarly we can prove the other two conclusions in
(ii).\end{proof}

\subsection{real $\Lambda_G(\sigma)-$representation}\label{realdef}

In this section  we discuss real $\Lambda_G(\sigma)-$\\
representation and its relation with the complex
$\Lambda_G(\sigma)-$representations introduced in Lemma \ref{cl}.
The main reference is \cite{BT} and \cite{FH}.

Let $G$ be a compact Lie group and $\sigma\in G^{n}_z$.
\begin{definition}A complex representation $\rho: G\longrightarrow Aut_{\mathbb{C}}(V)$ is said to be self dual if it is isomorphic to its complex dual
$\rho^*: G\longrightarrow Aut_{\mathbb{C}}(V^*)$ where
$V^*:=Hom_{\mathbb{C}}(V, \mathbb{C})$ and
$\rho^*(g)=\rho(g^{-1})^*$. \end{definition}

For any compact Lie group, we use $RO(G)$ to denote the real
representation ring of $G$. We have the real version of Lemma \ref{cl} below. The proof of Lemma \ref{clr} is left to the readers.
\begin{lemma}Let $\sigma\in G^{n}_z$. Then the map $\pi^*: RO(\mathbb{T}^n)\longrightarrow RO(\Lambda_G(\sigma))$
exhibits $RO(\Lambda_G(\sigma))$ as a free $RO(\mathbb{T}^n)-$module.

In particular there is an $RO(\mathbb{T})-$basis of
$RO(\Lambda_G(\sigma))$ given by irreducible real representations
$\{V_{\Lambda}\}$. There is a bijection between $\{V_{\Lambda}\}$
and the set $\{\lambda\}$ of irreducible real representations of
$C_G(\sigma)$. When $\sigma$ is trivial, $V_{\Lambda}$ has the
same underlying space $V$ as $\lambda$. When $\sigma$ is
nontrivial,
$V_{\Lambda}=((\lambda\otimes_{\mathbb{R}}\mathbb{C})\odot_{\mathbb{C}}(\eta_1\otimes\cdots\otimes \eta_n))\oplus
((\lambda\otimes_{\mathbb{R}}\mathbb{C})\odot_{\mathbb{C}}(\eta_1\otimes\cdots\otimes \eta_n))^*$
where each $\eta_i$ is a complex $\mathbb{R}-$representation such that
$(\lambda\otimes_{\mathbb{R}}\mathbb{C})(\sigma_i)$ acts on
$V\otimes_{\mathbb{R}}\mathbb{C}$ via the scalar multiplication by
$\eta_i(1)$. The dimension of $V_{\Lambda}$ is twice as that of
$\lambda$.\label{clr}\end{lemma}
As in (\ref{vgrepresentation}), we can construct a functor
$(-)^{\mathbb{R}}_{\sigma}$ from the category of real $G-$representations to the category of real
$\Lambda_G(\sigma)-$representations with
\begin{equation}(V)^{\mathbb{R}}_{\sigma} =
(V\otimes_{\mathbb{R}}\mathbb{C})_{\sigma}\oplus
(V\otimes_{\mathbb{R}}\mathbb{C})^{*}_{\sigma}.\label{realfaithv}\end{equation}
\begin{proposition}
Let $V$ be a faithful real $G-$representation. For each $\sigma\in
G^{n}_z$,
$(V)^{\mathbb{R}}_{\sigma}$ is a faithful real
$\Lambda_G(\sigma)-$representation. \label{farithreallambda}
\end{proposition}
\begin{proof}
Let $[a, t]\in \Lambda_G(\sigma)$ be an element acting trivially
on $(V)^{\mathbb{R}}_{\sigma}$. Assume $t\in [0, 1)$. Let $v\in
(V\otimes_{\mathbb{R}}\mathbb{C})_{\sigma}$ and let $v^*$ denote
its correspondence in
$(V\otimes_{\mathbb{R}}\mathbb{C})^*_{\sigma}$. Then $[a, t]\cdot
(v+v^*)=(a e^{2\pi i mt}+a e^{-2\pi i m t})(v+v^*)= v+v^*$ where
$m $ is a nonzero number determined by $\sigma$. Thus  $a$ is equal to
both $e^{2\pi i mt} I$, and $e^{-2\pi i mt}I$. Thus $t=0$ and $a$
is trivial.

So $(V)^{\mathbb{R}}_{\sigma}$ is a faithful real
$\Lambda_G(\sigma)-$representation.\end{proof}

\begin{proposition} Let $H$ and $G$ be two compact Lie groups. Let $\sigma \in G^{n}_z$ and $\tau\in H^{n}_z$. Let $V$ be a real $G-$representation and $W$ a real $H-$representation.

(i) We have the isomorphisms of representations $(V\oplus
W)^{\mathbb{R}}_{(\sigma, \tau)}=(V^{\mathbb{R}}_{\sigma}\oplus
                     W^{\mathbb{R}}_{\tau})$ as $\Lambda_{G\times H}(\sigma, \tau)\cong
                     \Lambda_G(\sigma)\times_{\mathbb{T}^n}\Lambda_H(\tau)-$representations.

(ii) Let $\phi: H\longrightarrow G$ be a group homomorphism. Let
$\phi_{\tau}: \Lambda_H(\tau)\longrightarrow
\Lambda_G(\phi(\tau))$ denote the group homomorphism obtained from
$\phi$. Then
$\phi_{\tau}^*(V)^{\mathbb{R}}_{\phi(\tau)}=(V)^{\mathbb{R}}_{\tau},$
as $\Lambda_H(\tau)-$representations.\end{proposition}
The proof
is left to the readers.


\begin{thebibliography}{10}




\bibitem{AHS}Matthew Ando, Michael J. Hopkins, and Neil P. Strickland: \textsl{Elliptic spectra, the Witten
genus and the theorem of the cube}. Invent. Math. , 146(3):595$\textendash$687,
2001.




\bibitem{BT}Theodor Br$\ddot{o}$cker, Tammo tom Dieck: \textsl{Representation of Compact Lie Groups},  Springer GTM \textbf{98} 1985.




\bibitem{Dev96}Jorge A. Devoto: \textsl{Equivariant elliptic homology and finite
groups}, Michigan Math. J. , 43(1):3$\textendash$32, 1996.




\bibitem{FH}William Fulton, Joe Harris: \textsl{Representation Theory, a first course},  Springer GTM \textbf{129} 1991.


\bibitem{Gan07}Nora Ganter: \textsl{Stringy power operations in Tate K-theory},
2007, available at arXiv: math/0701565.





\bibitem{HuanSpec} Zhen Huan: \textsl{Quasi-elliptic cohomology and its Spectrum},
available at arXiv:1703.06562.

\bibitem{Huanthesis}Zhen Huan: \textsl{Quasi-elliptic cohomology}, Thesis (Ph.D.)¨University of Illinois at Urbana-Champaign. 2017. 290 pp. \url{http://hdl.handle.net/2142/97268}.



\bibitem{Huansurvey} Zhen Huan: \textit{Quasi-elliptic cohomology I}. Advances in Mathematics. Volume 337, 15 October 2018, Pages 107-138.

\bibitem{Huanquasi} Zhen Huan: \textit{Quasi-theories}. available at arXiv:1809.06651.



\bibitem{KM85}Nicholas M. Katz and Barry Mazur: \textsl{Arithmetic moduli of elliptic curves}, Annals of Mathematics Studies, vol. 108, Princeton University Press, Princeton, NJ,
1985.



\bibitem{MM}M.A.Mandell, J.P.May: \textsl{Equivariant orthogonal spectra and
S-modules}, Mem.,Amer. Math. Soc. 159 (2002), no. 755, x+108 pp.







\bibitem{Rez11}Charles Rezk: \textsl{Quasi-elliptic cohomology}, unpublished manuscript, 2011.







\bibitem{SS}Stefan Schwede: \textsl{Global Homotopy Theory},
v0.23/April 30, 2015, Preliminary and incomplete version, \url{http://www.math.uni-bonn.de/people/schwede/global.pdf}.


\end{thebibliography}
\end{document}